\documentclass[a4paper, 10pt]{article}
\author{Femke Douma\\
Dept.~of Mathematical Sciences, Durham University,\\ Science Laboratories, South Road, Durham, DH1 3LE, UK\\ email: femke.douma@durham.ac.uk}
\title{A Lattice Point Problem on the Regular Tree}
\date{October 2009}
\usepackage[dvips]{graphicx}
\usepackage{amsmath}
\usepackage{amssymb}
\usepackage[all]{xy}
\usepackage{psfrag}
\newtheorem{thm}{Theorem}
\newtheorem{dfn}{Definition}
\newtheorem{lem}[dfn]{Lemma}

\newenvironment{pf}{\noindent \textsc{Proof} \rm}{\mbox{} \hfill $\square$}

\begin{document}

\maketitle

\begin{abstract}
Heinz Huber (1956) considered the following problem on the the hyperbolic plane $\mathbb{H}$. Consider a strictly hyperbolic subgroup of automorphisms on $\mathbb{H}$ with compact quotient, and choose a conjugacy class in this group. Count the number of vertices inside an increasing ball, which are images of a fixed point $x\in\mathbb{H}$ under automorphisms in the chosen conjugacy class, and describe the asymptotic behaviour of this number as the size of the ball goes to infinity. We use a well-known analogy between the hyperbolic plane and the regular tree to solve this problem on the regular tree.
\end{abstract}

\section{Introduction}\label{intro}

Let $X$ be a regular tree of degree $q+1\geq3$ with vertex set $V(X)$. Let $T:X\rightarrow X$ be a nontrivial hyperbolic automorphism or \emph{translation} on $X$, which is defined as an element in $\text{Aut}(X)$ with no fixed vertices or edges (see for example \cite[Chapter 1]{figa}). Denote by $d(x,Tx)$ the combinatorial \emph{distance} from a vertex $x \in V(X)$ to its image vertex under $T$. Let $\Gamma$ be a group of translations in $\text{Aut}(X)$ such that $G=X/\Gamma$ is a \emph{finite, simple, non-bipartite} $(q+1)$-regular graph, and let $\mathcal{K}$ be the conjugacy class of a nontrivial element $T_0\in\Gamma$. For $t\in\mathbb{R}$, $x\in V(X)$, we define
\begin{equation}\label{N}
N_{\mathcal{K}} (x,t) = \# \big\{ T \in \mathcal{K} : d(x,Tx) \leq t \big\}
\end{equation}
In this article we study the limiting behaviour of this counting function for increasing $t$. To do this we define a function $G_{\mathcal{K}}: V(X) \times \mathbb{C} \rightarrow \mathbb{C}$ for fixed $\mathcal{K}$ as follows:
\begin{equation}\label{G1}
G_{\mathcal{K}} (x,s) = \sum_{T \in \mathcal{K}} q^{-d(x,Tx)s}
\end{equation}
\begin{lem}\label{meromorphic}
The function $G_{\mathcal{K}} (x,s)$ as defined above \eqref{G1} is absolutely convergent for $\text{Re}(s)>2$.
\end{lem}
\begin{pf}
Note $\sum_{T \in \mathcal{K}} q^{-d(x,Tx)s} < \sum_{y\in V(X)} q^{-d(x,y)s}$. We can rewrite the last sum in terms of spheres $S_n(x)=\{y\in V(X) : d(x,y)=n\}$ as $\sum_{n=0}^\infty |S_n(x)| q^{-ns}$. Observe that $|S_n(x)|=(q+1)q^{n-1}$ for $n\geq1$ and use this to find $$\sum_{n=0}^\infty |S_n(x)| q^{-ns} \leq \sum_{n=0}^\infty 2q^nq^{-ns} = 2\sum_{n=0}^\infty q^{n(1-s)}$$ which converges for $\text{Re}(s)>2$ by the geometric series formula.
\end{pf}

In fact even more is true. In due course we will see that $G_\mathcal{K}$ has a meromorphic extension to $\mathbb{C}$, which is holomorphic for $\text{Re}(s)>\frac{1}{2}$.

The \emph{axis} $\mathfrak{a}(T)$ of a non-trivial translation $T$ on $V(X)$ is the unique geodesic which is mapped to itself by $T$. The \emph{displacement length} $\mu(T)$ of $T$ is given by
$$\mu(T) = \min_{x \in V(X)} d(x,Tx) \geq 1 $$
Note that $\mu(T) \in \mathbb{N}$. It is easy to see that the minimum is attained exactly for those vertices that lie in $\mathfrak{a}(T)$, as the vertices in $\mathfrak{a}(T)$ are shifted along the axis by $\mu(T)$ vertices under the action of $T$. Let $\delta(x,T)$ be the distance from a point $x$ to the axis of a translation $T$, that is
$$\delta(x,T) = \min_{y \in \mathfrak{a}(T)} d(x,y)$$
We then observe
\begin{equation}\label{distance}
d(x,Tx) = \mu(T) + 2 \delta(x,T)
\end{equation}

An element $P \in \Gamma$ is called \emph{primitive} if $\nexists \ Q \in \Gamma$ and $n>1$ such that $P=Q^n$. For every $T\in\Gamma$, we can find a unique $P \in \Gamma$ so that we can write $T=P^k$ for some $k\geq1$, and we call this the \emph{standard representation} of $T$. Now write $k=\nu(T)$ and call it the \emph{multiplicity} of $T$. Clearly $\mathfrak{a}(T) = \mathfrak{a}(P)$, and $\mu(T) = \mu(P) \cdot \nu(T)$. It is easy to prove that primitivity, $\mu(T)$ and $\nu(T)$ are invariant under conjugation in $\Gamma$.

We can now state our main theorem.

\begin{thm}\label{t1}
Let $N_{\mathcal{K}}(x,n)$ be defined as in equation \eqref{N} for positive integers $n$. Then as $n\rightarrow\infty$, where $n-\mu(\mathcal{K})$ is even,
$$N_{\mathcal{K}}(x,n) \sim q^{\frac{n-\mu(\mathcal{K})}{2}} \frac{ \mu(\mathcal{K})}{\nu(\mathcal{K}) |G|}$$
with $\mu(\mathcal{K})$ and $\nu(\mathcal{K})$ as defined above. By $N_{\mathcal{K}}(x,n)\sim Cq^{n/2}$ we mean that $N_{\mathcal{K}}(x,n)q^{-n/2}\rightarrow C$ as $n$ goes to infinity.
\end{thm}

We have to rule out the case where $n-\mu(\mathcal{K})$ is odd in the theorem because $N_\mathcal{K}(x,n)$ is a step function, which changes values only at points where $n-\mu(\mathcal{K})$ is even due to equation \eqref{distance}.

We are dealing with a type of problem that counts the number of lattice points on a graph inside an increasing ball. This type of problem is often solved using a discrete version of Selberg's trace formula on a regular tree (see \cite{selberg} and for the discrete version \cite{ahumada}, \cite{brooks}, \cite{t&w} or \cite{v&n}). In this particular case however, we are not able to use the trace formula due to the following. The trace formula is obtained using a $\Gamma$-invariant function $G(x,y)$ (a so-called point-pair invariant) such that
\begin{equation}\label{point-pair}
G(\gamma x,y)=G(x,\gamma y)=G(x,y)\ \forall\ \gamma\in\Gamma
\end{equation}
(see \cite{selberg} or \cite{t&w}). A natural choice would be $G_{\mathcal{K}}(x,y,s) = \sum_{T\in\mathcal{K}} q^{-d(x,Ty)s}$, but one easily checks that this function doesn't satisfy equation \eqref{point-pair} for arbitrary $x,y$. Instead, in this article we follow the method of Huber \cite{huber} using the function defined in equation \eqref{G1} which satisfies $G_{\mathcal{K}}(x,s) = G_{\mathcal{K}}(\gamma x,s)$ for all $\gamma\in\Gamma$.

Some other lattice point problems on the hyperbolic plane are easier to prove on the regular tree. Take for example the full lattice point problem, where we count $N_\Gamma(x,y,n) = \# \{ \gamma \in \Gamma : d(x,\gamma y)\leq n\}$. A consequence of the spherical average result in \cite{douma} for spheres $S_n(x_0)$ of vertices at distance exactly $n$ from $x_0$, is the asymptotic $\# \{ \gamma \in \Gamma : \gamma y \in S_n(x_0) \} \rightarrow \frac{|S_n(x_0)|}{|V(X/\Gamma)|}$. Using this and the fact that a ball of radius $n$ is the disjoint union of the spheres $\{S_r(x_0),\ 0\leq r \leq n\}$, we obtain
$$N_\Gamma(x,y,n)\rightarrow \frac{|B_n(x_0)|}{|V(X/\Gamma)|}$$
where we note that $|B_n(x_0)|$ is approximately equal to $\frac{q+1}{q-1}q^n$. For a proof of the corresponding problem on the hyperbolic plane, see \cite{patterson} or \cite[p 261]{buser}. Here we have
$$N_\Gamma(x,y,t) \rightarrow \frac{\pi}{\text{area}(M)}\ e^t$$
where $M$ is the manifold obtained by taking the quotient $\mathbb{H}/\Gamma$.

Another lattice point problem is that of counting the number of primitive closed paths of length at most $n$ in a finite regular graph $X/\Gamma$, which is a discrete analogue of counting primitive closed geodesics of length $\leq n$ on a closed hyperbolic surface. The result in the discrete case can be found in \cite[p 71]{sunada}, where we find that this number $\pi(n)=\# \{\text{prime cycles of length}\leq n \}$ has the following asymptotic behaviour
\begin{itemize}
\item when $X$ is non-bipartite, $\pi(n)\sim\frac{q^n}{n}$ as $n\rightarrow \infty$
\item when $X$ is bipartite, $\pi(n)\sim2\frac{q^n}{n}$ where $n$ is even and goes to $\infty$
\end{itemize}
The continuous analogue is the Prime Number Theorem for compact hyperbolic surfaces (see for example \cite[Theorem 9.4.14]{buser}). See also \cite{p&s} for a counting result of closed geodesics in negatively curved manifolds under homological constraints. A corresponding result for regular graphs can be found in \cite[p 72]{sunada}.

\section{Proof of the Theorem}

The discrete Laplacian is an operator which acts on functions on the vertices of any $(q+1)$-regular graph $G$ as follows:
\begin{equation}\label{laplacian}
\Delta f(x) = \frac{1}{q+1} \sum_{d(x,y)=1} f(y)
\end{equation}
Assume that the graph $G=X/\Gamma$ has $|V(G)|=N+1$ vertices. Then there is an orthonormal basis $\{\varphi_i\}_{i=0}^N$ of eigenfunctions of the Laplacian on the graph $G$. We use the fact that $\Delta = \frac{1}{q+1} A$ where $A$ is the adjacency matrix of $G$, to see that this basis has $N+1$ elements, and that we can order the eigenvalues such that $1=\lambda_0 > \lambda_1 \geq \ldots \geq \lambda_N > -1$ (note the eigenvalue $-1$ is excluded as $G$ is not bipartite: see for example \cite[Lemma 1.8]{chung}).

Using the canonical projection map $\pi:X\rightarrow G$, we can lift $\{\varphi_i\}_{i=0}^N$ to a set $\{ \tilde\varphi_i = \varphi_i\circ\pi \} _{i=0}^N$ of functions on $X$. It is easy to show that each $\tilde\varphi_i$ is an eigenfunction of the Laplacian on $X$ with the same eigenvalue $\lambda_i$. By definition the functions $\tilde\varphi_i$ are $\Gamma$-invariant. 

Recall we defined a function $G_\mathcal{K}(x,s)$ on $X$ which satisfies $G_\mathcal{K}(x,s) = G_\mathcal{K}(\gamma x,s)$ for all  $\gamma\in\Gamma$, so $G_\mathcal{K}$ is $\Gamma$-invariant. We can therefore view it as a function on the quotient $V(X/\Gamma)\times\mathbb{C}=V(G)\times\mathbb{C}$, where we call it $g_\mathcal{K}(x,s)$. Write $g_\mathcal{K}(x,s) = \sum_{i=0}^N F_i(s) \varphi_i(x)$. The Fourier coefficients $F_i(s)$ are given by $F_i(s) = \sum_{x\in V(G)} g_\mathcal{K}(x,s) \varphi_i(x)$. Now `lift' $g_\mathcal{K}(x,s)$ back up to $X$ and get
\begin{equation}\label{G_K}
G_\mathcal{K}(x,s) = \sum_{i=0}^N F_i(s)\tilde\varphi_i(x)
\end{equation}
where
\begin{equation}\label{fourier}
F_i(s) = \sum_{x\in\mathfrak{F}} G_\mathcal{K}(x,s) \tilde\varphi_i(x)
\end{equation}
for a fundamental domain $\mathfrak{F}\subset V(X)$ of $\Gamma$ on $X$. Note that there is a canonical one-to-one correspondence between $\mathfrak{F}$ and $V(G)$.

Choose and fix a translation $T^* \in \mathcal{K}$ with the standard representation $T^* = P^{\nu(\mathcal{K})}$ for some primitive $P$. Let \mbox{$H = \langle P \rangle$} be the subgroup of $\Gamma$ generated by $P$. We can write $\Gamma$ as a disjoint union of right cosets of $H$, i.e.
$$ \Gamma = \bigcup_{n=1}^{\infty} H A_n$$
with a fixed set $\{ A_n \}_{n=1}^\infty \subset \Gamma$. Then the elements $A_n^{-1} T^* A_n = T_n$ are pairwise disjoint, and run through all of $\mathcal{K}$ as $n=1,2,\ldots,\infty$. Define
$$\mathfrak{F}^* = \bigcup_{n=1}^\infty A_n(\mathfrak{F})$$
where $A_n(\mathfrak{F})=\{ A_n x : x\in\mathfrak{F} \}$. One easily checks that $\mathfrak{F}^*$ is a fundamental domain of the cyclic group $H$.

\begin{lem}\label{G*}
The Fourier coefficients $F_i(s)$ are given by
$$ F_i (s) = q^{-s\mu(\mathcal{K})} \sum_{x \in \mathfrak{F}^*} q^{-2s\delta(x,P)} \tilde \varphi_i (x)$$
\end{lem}
\begin{pf}
Let $k=\nu(\mathcal{K})$, and use equation \eqref{fourier} and the observations above to get
\begin{align}
F_i(s) & = \sum_{x\in \mathfrak{F}} G_{\mathcal{K}} (x,s) \tilde \varphi_i (x) = \sum_{x \in \mathfrak{F}} \sum_{n=1}^\infty q^{-d(x,A_n^{-1} P^{k} A_n x) s} \tilde \varphi_i (x) \nonumber \\
 & = \sum_{n=1}^\infty  \sum_{x \in A_n(\mathfrak{F})} q^{-d(x,P^k x)s} \tilde \varphi_i (x) = \sum_{x\in\mathfrak{F}^*} q^{-d(x,P^{k}x)s} \tilde\varphi_i(x) \nonumber \\
 & = q^{-sk\mu(P)} \sum_{x \in \mathfrak{F}^*} q^{-2s\delta(x,\mathfrak{a}(P))} \tilde \varphi_i (x) \label{e2}
\end{align}
where the last equality is due to equation \eqref{distance}. Note $k\mu(P) = \mu(\mathcal{K})$.
\end{pf}

Note that $\delta(x,P) = \delta(P^n x,P)$ for any integer $n$, as the axes of $P$ and $P^n$ coincide, and $\tilde\varphi_i(x) = \tilde\varphi_i(P^n x)$ as $P^n \in \Gamma$ and $\tilde\varphi_i(x)$ is $\Gamma$-invariant. This means the terms in the sum in equation \eqref{e2} is invariant under $H = \langle P \rangle$. Hence we can replace $\mathfrak{F}^*$ in the sum of \eqref{e2} by \emph{any} fundamental domain of $H$. Take a segment of $\mathfrak{a}(P)$ of length $\mu(P)$ and all branches emanating from the vertices in this segment, excluding the two branches that emenate from the vertices at the ends of the segment in the direction of the axis. The vertices in this set form a new fundamental domain for $H$, which we call $\mathfrak{F}_P$, and using the fact we can interchange fundamental domains of $H$ shown above, we now sum over the vertices in $\mathfrak{F}_P$ instead of $\mathfrak{F}^*$ in equation \eqref{e2}.

\begin{figure}[ht]
\centering
\psfrag{a(P)}{$\mathfrak{a}(P)$}
\psfrag{X}{$X$}
\psfrag{X/H}{$X/H$}
\includegraphics[height=4.5cm]{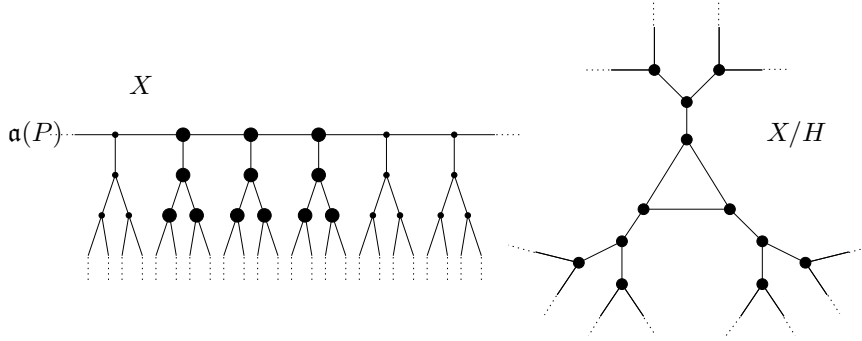}
\caption{The bold vertices in $X$ belong to $\mathfrak{F}_P$}
\end{figure}

The structure of the quotient graph $X/H$ can easily be deduced from $\mathfrak{F}_P$ (see for example Figure 1, where $\mu(P)=3$). We now want to transfer the functions $\tilde\varphi_i$ from $\mathfrak{F}_P$ to functions on the vertices of $X/H$, and to do this we use the obvious one-to-one correspondence between the vertex sets $\mathfrak{F}_P$ and $V(X/H)$. Note that the edge relations are preserved, and call the new function $\overline{\varphi_i}:V(X/H)\rightarrow\mathbb{R}$ for $i=0,\ldots,N$. It is easy to show these are eigenfunctions of the Laplacian on $X/H$. Using these definitions and equation \eqref{e2} we now have
\begin{align}
F_i(s) & = q^{-s\mu(\mathcal{K})} \sum_{x\in\mathfrak{F}_P} q^{-2s\delta(x,P)} \tilde\varphi_i(x) \nonumber \\
 & = q^{-s\mu(\mathcal{K})} \sum_{x\in V(X/H)} q^{-2s\delta'(x,P)} \overline\varphi_i(x) \label{e3}
\end{align}
where $\delta'(x,P)$ is the distance from the vertex $x$ to the central loop in $X/H$, which is exactly equal to $\delta(x,P)$ on $X$ (the central loop in Figure 1 is a triangle).

Define \emph{levels} in $X/H$ as follows:
$$L_n=L_n(X/H) = \{ x \in V(X/H) : \delta'(x,P)=n \} \quad \text{for} \ n\geq0$$
The \emph{radial average} of a function $f:V(X/H)\rightarrow\mathbb{R}$ with respect to these levels is $\frac{1}{|L_n|} \sum_{x \in L_n} f(x)$. A straightforward calculation shows that the radial average of $\overline{\varphi}_i (x)$ gives an eigenfunction $\Phi_i(x)$ of the Laplacian on $X/H$ with eigenvalue $\lambda_i$. Clearly $\Phi_i(x)=\Phi_i(y)$ whenever $\delta'(x,P)=\delta'(y,P)=n$, so we shall write $\Phi_i(n)$ from now on, where $n\in\mathbb{Z}^{\geq0}$. These $\Phi_i(n)$ are similar to the spherical functions with eigenvalue $\lambda_i$ on the tree as defined for example in \cite{t&w}. Use the facts that $V(X/H) = \bigcup_{n=0}^\infty L_n$ and $\sum_{x \in L_n} \overline \varphi_i (x) = |L_n| \Phi_i(n)$ in equation \eqref{e3} to obtain
\begin{equation}\label{e4}
F_i(s) = q^{-s\mu(\mathcal{K})} \sum_{n=0}^\infty |L_n| q^{-2sn} \Phi_i(n)
\end{equation}

\begin{lem}\label{lua}
There are constants $\alpha_i^\pm$, $u_i^\pm$ depending only on $\varphi_i$, $\lambda_i$ and $q$ such that
$$F_i(s) = q^{-s\mu(\mathcal{K})} \frac{\mu(\mathcal{K})}{\nu(\mathcal{K})} \frac{q-1}{q} \Big( \frac{u_i^+ + u_i^-}{q-1} + \frac{u_i^+}{1-\alpha_i^+ q^{1-2s}} + \frac{u_i^-}{1-\alpha_i^- q^{1-2s}} \Big)$$
except in the case that $\lambda_i = \pm \frac{2\sqrt{q}}{q+1}$, where we obtain
\begin{align}
F_i(s) = \frac{1}{q} \ q^{-s\mu(\mathcal{K})} & \frac{\mu(\mathcal{K})}{\nu(\mathcal{K})} \Phi_i(0) + q^{-s\mu(\mathcal{K})} \frac{\mu(\mathcal{K})}{\nu(\mathcal{K})} \frac{q-1}{q} \frac{\Phi_i(0)}{1\mp q^{1/2-2s}} \label{fourier D=0} \\
 & + q^{-s\mu(\mathcal{K})} \frac{\mu(\mathcal{K})}{\nu(\mathcal{K})} \frac{q-1}{q} \Phi_i(0) \big(\frac{\sqrt{q}-1}{\sqrt{q}+1}\big)^{\pm1}  \frac{\pm q^{1/2-2s}}{(1\mp q^{1/2-2s})^2} \nonumber
\end{align}
\end{lem}
\begin{pf}
Assume first that $\lambda_i\neq\pm\frac{2\sqrt{q}}{q+1}$. Note that $|L_0| = \mu(P)=\frac{\mu(\mathcal{K})}{\nu(\mathcal{K})}$, and $|L_n| = q^{n-1} (q-1) \frac{\mu(\mathcal{K})}{\nu(\mathcal{K})}$ for $n\geq1$. We use $\Delta \Phi_i(n) = \lambda_i \Phi_i(n)$ to find the recursion relation
\begin{equation}\label{recursion}
\Phi_i(n+1) = \frac{q+1}{q}\lambda_i\Phi_i(n) - \frac{1}{q}\Phi_i(n-1)\quad \text{ for } n\geq 1
\end{equation}
and $(q+1)\lambda_i\Phi_i(0) = (q-1)\Phi_i(1)+2\Phi_i(0)$. This leads to $\Phi_i(n) = u_i^+ (\alpha_i^+)^n + u_i^- (\alpha_i^-)^n$ for constants $\alpha_i^\pm$ and $u_i^\pm$ defined by
\begin{equation}\label{alpha}
\alpha_i^\pm = \frac{q+1}{2q}\lambda_i \pm \frac{\sqrt{(q+1)^2\lambda_i^2-4q}}{2q}
\end{equation}
\begin{equation}\label{u}
u_i^\pm = \bigg(\frac{1}{2} \pm \frac{(q+1)^2\lambda_i-4q}{2(q-1)\sqrt{(q+1)^2\lambda_i^2-4q}} \bigg) \cdot \Phi_i(0)
\end{equation}
Note that $\alpha_i^+\neq\alpha_i^-$ as the square root is nonzero due to the exclusion of $\lambda_i=\pm\frac{2\sqrt{q}}{q+1}$. Using the geometric series formula twice in equation \eqref{e4}, taking care with the $n=0$ term, we obtain
\begin{align}
F_i(s) = \frac{1}{q} \ q^{-s\mu(\mathcal{K})} \frac{\mu(\mathcal{K})}{\nu(\mathcal{K})} \Phi_i(0) & + q^{-s\mu(\mathcal{K})} \frac{\mu(\mathcal{K})}{\nu(\mathcal{K})} \frac{q-1}{q} \frac{u_i^+}{1-\alpha_i^+ q^{1-2s}} \nonumber \\
 & + q^{-s\mu(\mathcal{K})} \frac{\mu(\mathcal{K})}{\nu(\mathcal{K})} \frac{q-1}{q}  \frac{u_i^-}{1-\alpha_i^- q^{1-2s}} \label{fourier final}
\end{align}
Observe that for the infinite sums to converge, we need $|\alpha_i^\pm q^{1-2s}|<1$. For $i\neq0$ it is easy to check that $|\lambda_i|<1$ and $|\alpha_i^\pm| <1$ for $i\neq0$, and that there is a real number $\sigma_0<\frac{1}{2}$ so that the sums obtained from equation \eqref{e4} converge for $\text{Re}(s)>\sigma_0$. However for the eigenvalue $\mu_0=1$ we have $\alpha_0^+=1$ and the infinite sum for $F_0(s)$ will only converge for $\text{Re}(s)>\frac{1}{2}$.

The square root in equation \eqref{alpha} is zero iff we have $\lambda_i=\pm\frac{2\sqrt{q}}{q+1}$, which implies $\alpha_i^+=\alpha_i^-=\alpha_i=\pm\frac{1}{\sqrt{q}}$. In this case $\Phi_i(n)= ( 1 + n\big(\frac{\sqrt{q}-1}{\sqrt{q}+1}\big)^{\pm1} ) \Phi_i(0)\alpha_i^n$. Using the geometric series formula and the series $\sum_{i=1}^\infty ix^i = \frac{x}{(1-x)^2}$, we obtain the required expression from \eqref{e4}. For the convergence of the two infinite sums obtained here we require $|q^{1/2-2s}|<1$ which implies $\text{Re}(s)>\frac{1}{4}$, which is consistent with the general case above.
\end{pf}

Calculating $\alpha_0^-$ and $u_0^\pm$ using $\widetilde \varphi_0(x) = \frac{1}{\sqrt{|G|}} \ \forall \ x\in V(G)$ we obtain
\begin{equation}\label{G final}
G_{\mathcal{K}} (x,s) = \frac{\mu(\mathcal{K})}{\nu(\mathcal{K})|G|q^{s\mu(\mathcal{K})+1}} \big( 1+ \frac{q-1}{1-q^{1-2s}} \big) + \sum_{i=1}^{N} F_i(s) \widetilde\varphi_i(x)
\end{equation}
This is a meromorphic extension of $G_{\mathcal{K}} (x,s)$ defined in \eqref{G1} to the complex plane, which is holomorphic for $\text{Re}(s)>\frac{1}{2}$.

To finish the proof of our theorem, we use a refined version of the Tauberian theorem by Wiener-Ikehara from \cite{g&v} (see also \cite[chapter III Theorem 5.4]{korevaar}), which is a refinement of the Tauberian theorem in \cite{wiener}. This theorem requires a function $f(s)$ which converges for $\text{Re}(s)>1$ and has a simple pole at $s=1$. We choose $f(s) = G_{\mathcal{K}} (x,\frac{s}{2})$. The residue of $f(s)$ at $s=1$ is
\begin{equation}\label{A}
\text{Res}(f(s), 1) = \lim_{s\rightarrow1} (s-1)G_{\mathcal{K}} (x,\frac{s}{2}) =  \frac{\mu(\mathcal{K})(q-1)}{\nu(\mathcal{K})|G|q^{(\mu(\mathcal{K})/2)+1}}  \ \frac{1}{\ln q} :=A
\end{equation}
That means $f(s) = \frac{A}{s-1} + g(s)$ for some function $g(s)$ which is analytic for $\text{Re}(s)>1$. Check that $g(s)$ is analytic for $\text{Re}(s)\searrow1$ when $|\text{Im}(s)|<\frac{2\pi}{\ln q}$. Indeed, $g(s)$ has poles wherever $f(s)$ does, except we have removed the pole at $s=1$. As $f(s)$ has poles on the line $l=\{s:\text{Im}(s)=1\}$ at $s=1+ki\frac{2\pi m}{\ln q}$ for any $m\in\mathbb{Z}$, $g(s)$ has no poles on $l$ for $|\text{Im}(s)|<\frac{2\pi}{\ln q}$. Note also $A>0$.

Recall $N_{\mathcal{K}}(x,t) = \# \{ T \in \mathcal{K} : d(x,Tx) \leq t \}$. For fixed $x\in V(X)$, let
$$S(t) = N_{\mathcal{K}} (x,2t) = \# \{ T \in \mathcal{K} : d(x,Tx) \leq 2t \}$$
This is a non-decreasing step-function, which vanishes for $t<0$.

We now have all the ingredients we need to apply the theorem, which in our notation reads as follows.

\begin{thm}
Let $S(t)$ vanish for $t<0$, be non-decreasing, continuous from the right and such that
$$f(s) = \int_{0}^\infty q^{-s t} d S(t), \quad s = \sigma_1 + i \sigma_2$$
exists for $\text{Re}(s)=\sigma_1 >1$. Suppose that for some number $\rho>0$, there is a constant $A$ (necessarily $\geq0$) such that the analytic function
$$g(s) = f(s) - \frac{A}{s-1}, \quad s = \sigma_1 + i \sigma_2, \quad \sigma_1 >1$$
converges to a boundary function $g(1+i \sigma_2)$ in $L^1 (-\rho < \sigma_2 < \rho)$ as $\sigma_1 \searrow 1$. Let $\tau$ be the supremum of all possible numbers $\rho$. Then
$$ \frac{2\pi/\tau}{e^{2\pi/\tau}-1} A \leq \liminf_{t\rightarrow\infty} q^{-t} S(t) \leq \limsup_{t\rightarrow\infty} q^{-t} S(t) \leq \frac{2\pi/\tau}{1-e^{-2\pi/\tau}} A $$
\end{thm}

In our case $\tau = 2\pi/\ln q$. Using this and equation \eqref{A} we obtain
\begin{equation}\label{inf sup}
\frac{1}{q} \ \frac{\mu(\mathcal{K})}{\nu(\mathcal{K}) |G| q^{\mu(\mathcal{K})/2}}  \leq \liminf_{t\rightarrow\infty} \ q^{-t}S(t) \leq \limsup_{t\rightarrow\infty} \ q^{-t}S(t) \leq \frac{\mu(\mathcal{K})}{\nu(\mathcal{K})|G|q^{\mu(\mathcal{K})/2}}
\end{equation}
These estimates no longer depend on the choice of $x$.

Notice that as a consequence of \eqref{distance} when $\mu(\mathcal{K})$ is even, $S(t)$ will jump only at integer values of $t$ (the case where $\mu(\mathcal{K})$ is odd works similarly, except jumps occur only when $t+\frac{1}{2}$ is an integer). In this case, writing $m=[t]$ or equivalently
\begin{equation}\label{m}
t=m+\varepsilon\text{ with }m\in\mathbb{Z}\text{ and }\varepsilon\in [0,1)
\end{equation}
we have
\begin{equation}\label{t=m}
S(t)=S(m) \ \forall \ \varepsilon \in [0,1)
\end{equation}
Letting $a=\frac{\mu(\mathcal{K})}{\nu(\mathcal{K})|G|q^{\mu(\mathcal{K}) / 2}}$ we obtain the following estimates for the $\liminf$ and $\limsup$ respectively:
\begin{equation}\label{inf}
\frac{a}{q} \leq \liminf_{m\rightarrow\infty,m\in\mathbb{N}} q^{-(m+\varepsilon)} S(m) = q^{-\varepsilon} \liminf_{m\rightarrow\infty,m\in\mathbb{N}} q^{-m} S(m) \leq a
\end{equation}
\begin{equation}\label{sup}
\frac{a}{q} \leq \limsup_{m\rightarrow\infty,m\in\mathbb{N}} q^{-(m+\varepsilon)} S(m) = q^{-\varepsilon} \limsup_{m\rightarrow\infty,m\in\mathbb{N}} q^{-m} S(m) \leq a
\end{equation}
As this must hold for all $\varepsilon \in [0,1)$, we obtain
\begin{equation*}
\liminf_{m\rightarrow\infty,m\in\mathbb{N}} q^{-m} S(m) = \limsup_{m\rightarrow\infty,m\in\mathbb{N}} q^{-m} S(m) = \lim_{m\rightarrow\infty,m\in\mathbb{N}} q^{-m} S(m) = a
\end{equation*}
This means that for large integers $n=2m$ we have an approximation of $S(m)$ and hence $N_{\mathcal{K}}(x,n)$ as follows
$$S(m) \sim  q^{m-\frac{\mu(\mathcal{K})}{2}} \frac{\mu(\mathcal{K})}{\nu(\mathcal{K})|G|} \qquad \qquad N_{\mathcal{K}}(x,n) = S(\frac{n}{2}) \sim q^{\frac{n-\mu(\mathcal{K})}{2}} \frac{\mu(\mathcal{K})}{\nu(\mathcal{K}) |G|}$$
where from the definition of $m$ we require $n-\mu(\mathcal{K})$ to be even (because when $\mu(\mathcal{K})$ is odd, we take $m\in\frac{1}{2}+\mathbb{Z}$ in equation \eqref{m}). Using equation \eqref{t=m} it is clear that for any real $t$ such that $\mu(\mathcal{K})+2r\leq t < \mu(\mathcal{K})+2r+2$ for a non-negative integer $r$, we have $N_\mathcal{K}(x,t) = N_\mathcal{K}(x,\mu(\mathcal{K})+2r)$.

\textsc{Remark:} 
We now discuss why we required that $G$ was non-bipartite. Most of the proof above can be used to show a weaker result, but Theorem \ref{t1} will not hold for bipartite $G$. The method of proof works for the bipartite case up to Lemma \ref{lua}, where we have to consider the effects on $G_\mathcal{K}(x,s)$ of $\lambda_N=-1$, the eigenvalue of the Laplacian which occurs exactly when $G$ is bipartite. The series for $F_N(s)$ requires $\text{Re}(s)>\frac{1}{2}$ to converge, and we obtain
\begin{align*}
G_{\mathcal{K}} (x,s) = & \frac{\mu(\mathcal{K})}{\nu(\mathcal{K})|G|q^{s\mu(\mathcal{K})+1}} \big( 1+ \frac{q-1}{1-q^{1-2s}} \big) + \sum_{i=1}^{N-1} F_i(s) \widetilde\varphi_i(x) \\
& + \frac{\mu(\mathcal{K})}{\nu(\mathcal{K})|G|q^{s\mu(\mathcal{K})+1}} \big(1+\frac{q-1}{1+q^{1-2s}}\big)
\end{align*}
As before the residue of $f(s)=G_\mathcal{K}(x,\frac{s}{2})$ at $s=1$ equals $A$, but now the function $g(s)=f(s)-\frac{A}{s-1}$ only converges to a boundary function for $|\text{Im}(s)|<\frac{\pi}{\ln q}=\tau$. The Tauberian theorem can still be applied, but only shows
\begin{align*}
\frac{1}{q^2} \frac{2q}{q+1} \ \frac{\mu(\mathcal{K})}{\nu(\mathcal{K})|G|q^{\mu(\mathcal{K})/2}}  \leq & \liminf_{t\rightarrow\infty} \ q^{-t}S(t) \\
\leq & \limsup_{t\rightarrow\infty} \ q^{-t}S(t) \leq \frac{2q}{q+1} \frac{\mu(\mathcal{K})}{\nu(\mathcal{K})|G|q^{\mu(\mathcal{K})/2}}
\end{align*}
The difference between this estimate and that in equation \eqref{inf sup} is that the left and right estimates differ by a factor $q^2$ instead of the factor $q$ obtained in the non-bipartite case, and it is no longer possible to deduce a precise limit from the estimates (note that $S(t)$ can jump at any integer value $t$). All we can say here is
\begin{align*}
\frac{1}{q^2} \frac{2q}{q+1} \ \frac{\mu(\mathcal{K})}{\nu(\mathcal{K})|G|q^{\mu(\mathcal{K})/2}}  \leq & \liminf_{t\rightarrow\infty} \ q^{-t/2}N_\mathcal{K}(x,t) \\
\leq & \limsup_{t\rightarrow\infty} \ q^{-t/2}N_\mathcal{K}(x,t) \leq \frac{2q}{q+1} \frac{\mu(\mathcal{K})}{\nu(\mathcal{K})|G|q^{\mu(\mathcal{K})/2}}
\end{align*}

Acknowledgements: The author wishes to thank N.~Peyerimhoff for many helpful discussions. This work forms part of the author's PhD research, which is supported by the EPSRC.

\end{document}